\documentclass[12pt]{amsart}
\usepackage{amsmath,amsthm,amsfonts,amscd,amssymb}

\theoremstyle{plain}

\newtheorem{thm}{Theorem}[section]
\newtheorem{defn}[thm]{Definition}
\newtheorem{lem}[thm]{Lemma}
\newtheorem{cor}[thm]{Corollary}

\newtheorem*{thm*}{Theorem}
\newtheorem*{rem*}{Remark}
\newtheorem*{lem*}{Lemma}

\newtheorem*{cor*}{Corollary}
\newtheorem*{prop*}{Proposition}

\newcommand{\ra}{\rightarrow}

\newcommand{\dra}{\dashrightarrow}

\newcommand{\ov}{\overline}

\newcommand{\bZ}{\mathbb{Z}}

\newcommand{\bP}{\mathbb{P}}

\newcommand{\bd}{\begin{equation*} \begin{diagram}}
\newcommand{\ed}{\end{diagram} \end{equation*}}

\begin{document}

\title{Motives of unitary and orthogonal homogeneous varieties}

%\thanks{This material is based upon work supported by the National
%Science Foundation under agreement No. DMS-0111298. Any opinions,
%findings and conclusions or reccomendations expressed in this material
%are those of the author and do not necessarily reflect the views of
%the National Science Foundation.}

\author{Daniel Krashen}

\begin{abstract}

Grothendieck-Chow motives of quadric hypersurfaces have provided many
insights into the theory of quadratic forms. Subsequently, the
landscape of motives of more general projective homogeneous varieties
has begun to emerge. In particular, there have been many results which
relate the motive of a one homogeneous variety to motives of other
simpler or smaller ones (see for example \cite{Kar:cell, CGM, Bro:BB,
NSZ:F4, SZ:Hom, CPSV, Zan:GSBV}. In this paper, we exhibit a
relationship between motives of two homogeneous varieties by producing
a natural rational map between them which becomes a projective bundle
morphism after being resolved. This allows one to use formulas for
projective bundles and blowing up (as in \cite{Man:Mot}) to relate the
motives of the two varieties. We believe that in the future this ideas
could be used to discover more relationships between other types of
homogeneous varieties.

\end{abstract}

\maketitle

\section{Introduction}

Let $F$ be a field of characteristic not $2$, and let $L/F$ be a
quadratic extension with Galois group $G = \left< \sigma \right>$. Let
$W$ be an $L$-vector space and $h$ an $L/F$ Hermitian form. We will
abbreviate this by saying that $(W, h)$ is an $L/F$ Hermitian space.

Both quadratic and Hermitian spaces give rise to projective
homogeneous varieties as follows. For a quadratic space $(W, q)$, we
let $V(q)$ denote the associated quadric hypersurface in $\bP(W)$
defined by the vanishing of $q$. This is a homogeneous space for the
orthogonal group $O(W, q)$.

If $(W,h)$ is an $L/F$-Hermitian space we denote by $V(h)$ the variety
of isotropic $L$-lines in $W$. This is an $F$-variety and may be
regarded as a closed subvariety of $tr_{L/F} \bP_L(W)$ which in turn is
a closed subvariety of the Grassmannian $Gr(2, W)$ of $2$-dimensional
$F$-subspaces of $W$. It is a homogeneous space for the unitary group
$U(W, h)$, and is a twisted form of the variety of flags consisting of
a dimension 1 and codimension 1 linear subspace in a $n$ dimensional
vector space (where $dim_F(W) = 2n$).

For a regular $F$-scheme $X$, we let $M(X)$ denote the
Grothendieck-Chow motive of $X$ as described in \cite{Man:Mot}. For a
motive $M$, we let $M(i)$ denote the $i$'th Tate twist of $M$. Our
main result is the following:

\begin{thm*}\textnormal{\textbf{\ref{main_thm}.}}
Let $(W,h)$ be an $L/F$-Hermitian space, and let $q_h$ be the
associated quadratic form on $W$. Then we have an isomorphism of
motives
\begin{equation*}
M(V(q_h)) \oplus \bigoplus_{i = 1}^{n-2} M(\bP_L(N))(i) \cong
M(V(h)) \oplus M(V(h))(1)
\end{equation*}
where $dim_F(N) = 2n$.
\end{thm*}
Comparing the dimension $0$ Chow groups of each side,
we obtain:

\begin{cor*}\textnormal{\textbf{\ref{main_cor}.}}
\begin{equation*}
A_0(V(h)) \cong A_0(V(q_h)) = \left\{
\begin{matrix}
\bZ && \text{if $q$ is isotropic} \\
2\bZ && \text{if $q$ is anisotropic,}
\end{matrix}
\right.
\end{equation*}
in particular, the degree map $A_0(V_h) \ra \bZ$ is injective.
\end{cor*}

\section{Projective Spaces}

\begin{defn}
Let $L$ be an \'etale $F$-algebra, and let $N$ be an $L$-module of
rank $m < \infty$. We define $\bP(N, L)$ to be $tr_{L/F} \bP_L(N)$
which may be regarded as a closed subvariety of $Gr(dim(L), N)$, where
$N$ is considered here as an $F$-space.
\end{defn}

\begin{defn}
$S(N) = \{[m] \in \bP_F(N) | ann_L(m) \neq 0\}$.
\end{defn}

We note the following easy lemma:

\begin{lem}
Let $L$ be a commutative \'etale $F$-algebra, and let $N$ be a finite
dimensional free $L$-module. Then we have a rational morphism $f = f_N
: \bP_F(N) \dra \bP(N, L)$ defined by taking a $1$-dimensional linear
subspace $l \subset N$ to its (usually $dim(L)$-dimensional) $L$-span,
$Ll$. This is well defined on the open complement of $S(N)$ in
$\bP_F(N)$.
\end{lem}

This result may be proved by descent. We omit the proof, but note that
it is essentially contained in the proof of the following result.

\begin{lem} \label{proj_case}
Suppose $L$ is either a quadratic \'etale extension of $F$. Let $N$ be
a finite rank free $L$-module, and let $S(N) = \{[m] \in
\bP_F(N) | ann_L(m) \neq 0\}$. Then there is a well defined morphism
$g_N : Bl_{S(N)} \bP_{F}(N) \ra \bP_L(N)$ extending the rational map
$f_N$. Furthermore, $g_N$ is a projective bundle of rank $1$.
\end{lem}

\begin{proof}
By descent, it suffices to consider the case where $F$ is
algebraically closed. In this case, we have $L = F e_1 \oplus F e_2$
and so $N = N_1 \oplus N_2$, where each $N_i = N e_i$ is a vector
space of finite dimension $n$. This gives
$$\bP_L(N) = \bP_F(N_1) \times \bP_F(N_2),$$ and
$$S(N) = \{[n_1; n_2] \in \bP_F(N) | \text{ either $n_i$ or $n_j$ is
$0$}\},$$ where the rational morphisms $\bP_F(N) \ra \bP_F(N_i)$ are
given by projections. One may check directly that the desired morphism
$g_N$ is just the resolution of the rational morphism. More precisely,
we have an isomorphism
\begin{multline*}
Bl_{S(N)} \bP_F(N) = \\ 
\left\{([n_1; n_2], [m_1], [m_2]) \in \bP_F(N) \times \bP_F(N_1)
\times \bP_F(N_2) \left| \
\begin{matrix}
n_i \text{ is in } \\
\text{the line } [m_i]
\end{matrix} \ 
\right.\right\}.
\end{multline*}
In particular, this is a $\bP^1$ bundle as the inverse image of
$[n_1], [n_2]$ is given by $\{[\lambda n_1 ; \mu n_2] | \lambda, \mu
\in F\}$, which is isomorphic as a variety with a Galois action to
$\bP_F(L)$.
\end{proof}

\begin{lem}
Let $L$ be a quadratic \'etale extension of $F$. Using the notation of
lemma \ref{proj_case}, there is an isomorphism $S(N) \cong \bP_L(N)$.
\end{lem}
\begin{proof}
Choose an $L$ basis $n_1, \ldots n_k$ for $N$. Since $L$ is \'etale,
we may write $L_{\ov{F}} = \ov{F}e_1 \oplus \ov{F}e_2$, where $\ov{F}$
is the separable closure of $F$. Then $S(N)_{\ov{F}}$ is the disjoint
union of the schemes $P_i$, $i = 1, 2$, where $P_i =
\bP_{\ov{F}}(W_i)$ and $W_i$ is the $\ov{F}$-span of the vectors $n_j
e_i$ as $j$ varies. Examining the $G_F$-Galois action, we see that we
have an isomorphism of $\ov{F}$-schemes between $S(N)_{\ov{F}}$ and
$\bP_{\ov{F}}(W) \times_{\ov{F}} Spec(L_{\ov{F}})$, where $W$ is a
$k$-dimensional $\ov{F}$ vector space. By descent, this gives $S(N) =
\bP_L(V)$, where $V$ is a rank $k$ free $L$-module, and of course we
may choose $V = N$.
\end{proof}

As a result of this we obtain information on the motive
$M(\bP_L(N))$:

\begin{thm}
Let $L/F$ be a quadratic \'etale extension. Then we have an isomorphism
in the category of Grothendieck-Chow motives:
\begin{multline*}
M(\bP_F(N)) \oplus \bigoplus_{i = 1}^{n-1} M(\bP_L(N))(i) \cong
M(\bP(N, L)) \oplus M(\bP(N, L))(1) \\ \cong
M(tr_{L/F} \bP(N)) \oplus M(tr_{L/F} \bP(N))(1)
\end{multline*}
where $dim_F(N) = 2n$.
\end{thm}
\begin{proof}
This follows from the formulas in \cite{Man:Mot} for the motive of a
projective bundle and for the motive of a blow up.
\end{proof}

\section{Quadratic and Hermitian Spaces}

Let $(W, h)$ be an $L/F$ Hermitian space. Write $h$ in diagonal form
as $h = \left<a_1, \ldots, a_n\right>$, taken with respect to a basis
$w_1, \ldots, w_n$ for $W$ over $L$. If $L = F(\beta)$ with $\beta^2 =
b \in F$, we may consider the basis $u_i = w_i, v_i = \beta w_i$ for
$W$ as an $F$-space. Note that we may define a quadratic form $q_h$ on
$W$ as an $F$-space by $v \mapsto q_h(v) = h(v, v)$. In this basis
this is exactly the form $\left<1, -b\right> \otimes \left<a_1,
\ldots, a_n\right>$.

We note the following observation:

\begin{lem}
Let $(W, h)$ be an $L/F$ Hermitian space. Then $w \in W$ is isotropic
for $h$ if and only if it is isotropic for $q_h$.
\end{lem}

Consequently, if we consider the map from the previous section:
$$f = f_W : \bP_F(W) \dra \bP(W, L),$$
we obtain by restriction a rational map
$$\phi : V(q_h) \dra V(h).$$

\begin{lem} \label{main_case}
Suppose $L$ is either a quadratic \'etale extension of $F$. Let $W$ be
a finite dimensional free $L$-module, and let 
$$S(W) = \{[m] \in \bP_F(W) | ann_L(m) \neq 0\}.$$ Then $S(W) \subset
V(q_h)$ and the morphism $g_W : Bl_{S(W)} \bP_{F}(W) \ra \bP_L(W)$
from lemma \ref{proj_case} restricts to a morphism $g_h : Bl_{S(W)}
V(q_h) \ra V(h).$
\end{lem}

\begin{proof}
Suppose that $S(W) \subset V(q_h)$. Then by the standard properties of
the blow up (\cite{Ful:IT}, B.6.10), we get a natural (inclusion)
morphism $Bl_{S(W)} V(q_h) \ra Bl_{S(W)} V(h)$, and so by restriction
a morphism $Bl_{S(W)} V(q_h) \ra \bP_L(W)$. Since this map agrees with
$\phi$ on a dense open set, it follows that the image is exactly $V(h)
\subset \bP_L(W)$.

It remains to show that $S(W) \subset V(q_h)$. To show this, it
suffices to assume that $\ov{F} = F$, and so $L \cong F \oplus F$, say
with idempotent basis $e_1, e_2$. In this case we may suppose that the
Galois involution on $L$ is $\sigma$ defined by $\sigma(e_1) = e_2,
\sigma(e_2) = e_1$, and we may assume $h = <1, \ldots, 1>$ with
respect to a given basis $w_i$. Choosing the basis $u_i = e_1 w_i, v_i
= e_2 w_i$ for $W$ as an $F$-space, we see that for a vector 
$$w = \sum_{i = 1}^n \lambda_i u_i + \mu_i v_i,$$
$w \in S(W)$ if and only if either $\lambda_i = 0$ all $i$ or $\mu_i =
0$ all $i$. Since $$q_h(w) = \sum \lambda_i \mu_i,$$ it immediately
follows that $S(W) \subset V(q_h)$ as desired.
\end{proof}

Since the morphism $g_h$ is obtained by restriction of $g_W$, it is
also a projective bundle morphism of relative dimension $1$, and again
we have from the motivic blow-up and projective bundle formulas
(\cite{Man:Mot}):

\begin{thm} \label{main_thm}
Let $(W,h)$ be an $L/F$-Hermitian space, and let $q_h$ be the
associated quadratic form on $W$. Then we have an isomorphism of
motives
\begin{equation*}
M(V(q_h)) \oplus \bigoplus_{i = 1}^{n-2} M(\bP_L(W))(i) \cong
M(V(h)) \oplus M(V(h))(1)
\end{equation*}
where $dim_F(W) = 2n$.
\end{thm}

In particular, comparing the dimension $0$ Chow groups of each side,
we obtain:

\begin{cor} \label{main_cor}
\begin{equation*}
A_0(V(h)) \cong A_0(V(q_h)) = \left\{
\begin{matrix}
\bZ && \text{if $q$ is isotropic} \\
2\bZ && \text{if $q$ is anisotropic,}
\end{matrix}
\right.
\end{equation*}
in particular, the degree map $A_0(V_h) \ra \bZ$ is injective.
\end{cor}

\bibliographystyle{alpha}
\bibliography{citations}

\end{document}